\documentstyle[12pt]{article}
\begin{document}

\title{Unions of Cockcroft two-complexes}
\author{W. A. Bogley \\ Oregon State University \\ Corvallis, OR  97331 USA}
\date{September, 1991\\ revised October, 1992}
\maketitle

\begin{abstract} A combinatorial group-theoretic hypothesis is presented that
serves as a necessary and sufficient condition for a union of connected
Cockcroft two-complexes to be Cockcroft.  This hypothesis has a component that
can be expressed in terms of the second homology of groups.  The hypothesis is
applied to the study of the third homology of groups given by generators and
relators.

\vspace{3mm}

\end{abstract}  

\section{Introduction}

It was observed by Cockcroft \cite{c} that if a two-complex $K$ (that
is, a two-dimensional CW complex) is a subcomplex of an aspherical
two-complex, then for any basepoint $z \in K$, the Hurewicz map \(h:\pi_2(K,z)
\rightarrow H_2K \) is trivial.  In other words, each spherical map \(S^2
\rightarrow K\)  is homologically trivial.  Two-complexes with this latter
property are therefore said to be {\em Cockcroft}.  In this note, necessary
and sufficient conditions are given for a union $K$ of subcomplexes $K_r$ and
$K_s$ to be Cockcroft.  A necessary condition is that each of $K_r$ and $K_s$
be Cockcroft.  An additional condtion is group-theoretic; it is of interest in
connection with a certain ``relative Hopf formula'' for third homology that
was introduced in \cite{bg}.

This introductory section discusses refinements and the significance of the
Cockcroft property.  The main result is proved in Section 2.  Examples and
applications to homology calculations are presented in Section 3.

The Cockcroft property has group-theoretic content.  Suppose that $K$ is
modeled on a presentation $\cal P = ({\bf x} : {\bf u} )$ for a group $G$.
For $u \in {\bf u}$, write $u = q^e$ where $e \geq 1$ and $q$ is not a proper
power.  If $K$ is Cockcroft, then $q$ determines an element of order exactly
$e$ in $G \cong \pi_1 K$\cite{hu,sjp91}.  Next, there is the exact Hopf
sequence \cite{ho1,ho2}

\[ 0 \rightarrow H_3 G \rightarrow {\bf Z}
\otimes_{{\bf Z}G} \pi_2 K \rightarrow H_2 K \rightarrow H_2 G \rightarrow 0\]

\noindent where the middle map is naturally induced by $h$.  It follows that
if $K$ is Cockcroft, then \(H_2G \cong H_2K\) is free abelian and $H_3 G \cong
{\bf Z} \otimes_{{\bf Z}G} \pi_2 K$.  From a combinatorial perspective, if $K$
is finite and Cockcroft, then $\cal P$ is {\em efficient}, in the sense that
the number of generators minus the number of relators is equal to the
difference of the torsion-free rank of $H_1 G$ minus the minimum number of
generators for $H_2 G$ \cite{d2}.  In particular, $K$ has minimum Euler
characteristic among all finite two-complexes with fundamental group
isomorphic to $G$ \cite{bd}.  Finally, the Cockcroft property has been used to
produce lower bounds for isoperimetric functions of group presentations
\cite{g,hr}.

Increasingly delicate versions of the Cockcroft property arise upon passage to
coverings of  $K$.  For a subgroup  $H \leq \pi_1 K$,  the two-complex  $K$
is  {\em $H-$Cockcroft}  if the lifted Hurewicz map $\bar{h} :
\pi_2K \rightarrow H_2\bar{K}$ is trivial, where $\bar{K} \rightarrow K$ is
the covering corresponding to $H$.  This is the same as saying that $\bar{K}$
is Cockcroft.  In particular, $K$ is aspherical if and only if $K$ is
$\{1\}$-Cockcroft.  For $H,H' \leq \pi_1K$, if $K$ is $H-$Cockcroft and some
$\pi_1K$-conjugate of $H$ is contained in $H'$, then $K$ is $H'$-Cockcroft.
(This is because conjugate subgroups arise by changing basepoints in covering
complexes.)

As an example, suppose that $K$ is modeled on a one-relator presentation
$({\bf x}: r)$, where $r \not= 1$ in the free group $F$ with basis ${\bf x}$.
Write $r = q^e$ where $q$ is not a proper power in $F$.  It follows from
Lyndon's Simple Identity Theorem \cite{l} that for a subgroup $H \leq \pi_1
K$, $K$ is $H$-Cockcroft if and only if $H$ contains each $\pi_1 K$-conjugate
of $q$.  The point is that the Identity Theorem can be interpreted as a
description of a generating set for $\pi_2 K$.  Generalizations of this
observation appear in \cite{gh1}.

These refined versions of the Cockcroft property first appeared in \cite{bd},
and have since received considerable attention
\cite{bds,d1,d2,h,gh1,gh2,sjp91,sjp92}.  They provide an intrinsic connection
between the subgroup structure of $\pi_1 K$ and the module structure of $\pi_2
K$.  It is shown in \cite{gh1} and in \cite{h} that if $K$ is Cockcroft, then
$\pi_1 K$ contains a minimal subgroup $H$ such that $K$ is $H$-Cockcroft.
These so-called {\em threshold} subgroups are studied extensively in
\cite{gh1,gh2}.  Examples of two-complexes with nonunique thresholds are
presented in \cite{sjp92}.

\section{The main result}

Suppose that a connected two-complex $K$ is given as the union of connected
subcomplexes $K_r$ and $K_s$, where $K_r \cap K_s = K^{(1)}$.  Let $F$ =
$\pi_1 K^{(1)}$ and let $R = \ker(F \rightarrow \pi_1K_r)$ and $S = \ker(F
\rightarrow \pi_1K_s)$ be the kernels of the inclusion-induced maps.  The
subgroups of $\pi_1 K = F/RS$ are all of the form $N/RS$, where $RS \subseteq
N \leq F$.  If $A$ and $B$ are subgroups of $F$, then $[A,B]$ denotes the
subgroup of $F$ generated by all commutators $[a,b] = aba^{-1}b^{-1}$, where
$a \in A$ and $b \in B$.

\vspace{3mm}
\noindent  {\bf Theorem} {\it For a subgroup  $N \leq F$  with  $RS \subseteq N$,  $K$  is  $N/RS$-Cockcroft if and only if
	\begin{enumerate}
		\item $K_r$ is $N/R$-Cockcroft,
		\item $K_s$ is $N/S$-Cockcroft, and
		\item $R \cap S \subseteq [R,N] \cap [S,N]$.
	\end{enumerate} }
\vspace{3mm}

\noindent  {\bf Proof:}  Let  $p:\bar{K} \rightarrow K$  be the covering of
$K$ corresponding to $N/RS$.  Let $p_r :\bar{K}_r \rightarrow K_r$ be the
restriction of $p$ to $\bar{K}_r = p^{-1}(K_r)$; $p_r$ is the covering of
$K_r$ corresponding to $N/R$.  If $K$ is $N/RS$-Cockcroft, then $\bar{K}$ is
Cockcroft, and so the subcomplex $\bar{K}_r$ is Cockcroft.  Similarly,
$\bar{K}_s = p^{-1}(K_s)$ is Cockcroft.  We may thus assume throughout that
$K_r$ is $N/R$-Cockcroft and that $K_s$ is $N/S$-Cockcroft.  The following
commutative diagram has exact rows and columns.  The top row is from
\cite[Theorem 1]{gr} (see also \cite[Corollary 3.2]{b} and \cite[Corollary
3.4]{bg}) and the middle row is from the Mayer-Vietoris homology sequence for
$\bar{K} = \bar{K}_r \cup \bar{K}_s$.  Exactness of the first two columns is
due to Hopf \cite{ho1,ho2}.

\vspace{7mm}

\begin{picture}(350,150)

%  This commutative diagram requires nine columns, having unit lengths
%	30	30	80	30	40	30	40	30	30.
%  It has five rows, with unit lengths  (from top to bottom)
%	15	50	15	50	15.
%  Horizontal arrow boxes are 30 x 15.  Text rows have height 15.

%  First text row

	\put(60,130){\makebox(80,15){$\pi_2K_r \oplus \pi_2\,K_s$}}
	\put(140,130){\makebox(30,15){$\rightarrow$}}
	\put(170,130){\makebox(40,15){$\pi_2\,K$}}
	\put(210,130){\makebox(30,17)[t]{$\stackrel{\eta}{\rightarrow}$}}
	\put(240,130){\makebox(30,15){$\frac{R \cap S}{[R,S]}$}}
	\put(270,130){\makebox(30,15){$\rightarrow$}}
	\put(300,130){\makebox(30,15){$0$}}

%vertical arrow row

	\put(100,120){\vector(0,-1){30}}
	\put(105,105){$0$}
	\put(190,120){\vector(0,-1){30}}
	\put(195,105){$\bar{h}$}
	\put(255,120){\vector(0,-1){30}}
	\put(260,105){$\nu^N$}

%  Second text row

	\put(0,65){\makebox(30,15){$0$}}
	\put(30,65){\makebox(30,15){$\rightarrow$}}
	\put(60,65){\makebox(80,15){$H_2\bar{K}_r \oplus H_2\bar{K}_s$}}
	\put(140,65){\makebox(30,15){$\rightarrow$}}
	\put(170,65){\makebox(40,15){$H_2\bar{K}$}}
	\put(210,65){\makebox(30,15){$\rightarrow$}}
	\put(240,65){\makebox(30,15){$H_1\bar{K}^{(1)}$}}

%vertical arrow row

	\put(100,55){\vector(0,-1){30}}
	\put(105,40){$\cong$}
	\put(190,55){\vector(0,-1){30}}
	\put(195,40){onto}

%  Third text row

	\put(60,0){\makebox(80,15){$H_2\,N/R \oplus H_2\, N/S$}}
	\put(140,0){\makebox(30,15){$\rightarrow$}}
	\put(170,0){\makebox(40,15){$H_2\,N/RS$}}

\end{picture}

\vspace{7mm}

As in \cite[Lemma 3.1]{bg}, the homomorphism $\nu^N$  can be identified with
the map 
	\[ \frac{R \cap S}{[R,S]} \rightarrow \frac{N}{[N,N]} \]
that is induced by the inclusion of  $R \cap S$  in  $N$.  Applying the snake
lemma, Hopf's formula for the second homology of groups shows that the bottom
row of the diagram can be extended to the left and rewritten as the exact
sequence
	\[ \ker \nu^N \stackrel{i}{\rightarrow} \frac{R \cap [N,N]}{[R,N]}
\oplus \frac{S \cap [N,N]}{[S,N]} \rightarrow \frac{RS \cap [N,N]}{[RS,N]}, \]
where
	\[ i(w[R,S]) = (w[R,N],w^{-1}[S,N]) \]

\noindent for each   $w \in R \cap S \cap [N,N]$.  To see this, a
description of the map $\eta$ is given in \cite[p. 49]{gr} and in
\cite[Remarks 2.6 and 3.6]{bg}.  Namely, given $w \in R \cap S$, there exist
singular discs $\alpha : (B^2 ,S^1 ) \rightarrow (\bar{K}_r ,\bar{K}^{(1)})$
and $\beta: (B^2 ,S^1 ) \rightarrow (\bar{K}_s ,\bar{K}^{(1)})$ with
boundaries $\alpha |_{S^1} = \beta |_{S^1}$ representing $w \in \pi_1
\bar{K}^{(1)} \cong N$.  One then has that

\[ [\alpha][\beta]^{-1} \in
\ker (\pi_2 (\bar{K},\bar{K}^{(1)}) \rightarrow \pi_1 \bar{K}^{(1)}) \cong \pi_2 K\]

\noindent and that $\eta([\alpha][\beta]^{-1}) = w[R,S]$.  If $w \in R \cap S
\cap [N,N]$ then the image of $[\alpha]$ under the relative Hurewicz map $\pi_2
(\bar{K}_r ,\bar{K}^{(1)}) \rightarrow H_2 (\bar{K}_r ,\bar{K}^{(1)})$ is an
element of $\ker (H_2 (\bar{K}_r ,\bar{K}^{(1)}) \rightarrow H_1
\bar{K}^{(1)}) \cong H_2 \bar{K}_r$.  Further, this element projects to
$w[R,N]$ under the Hopf map $H_2 \bar{K}_r \rightarrow (R \cap [N,N])/[R,N]$.
Similar remarks apply to $[\beta]$ and the description of $i$ follows
directly.

A diagram chase now reveals that
	\[ K\  {\rm is}\  N/RS{\rm -Cockcroft}\  \Leftrightarrow \bar{h} = 0
\Leftrightarrow \nu^N = 0 \ {\rm and}\  i = 0. \]
The definitions of  $\nu^N$  and  $i$  thus imply that  $K$  is
$N/RS$-Cockcroft if and only if
	\[ R \cap S \subseteq [N,N] \  {\rm and}\   R \cap S \cap [N,N]
\subseteq [R,N] \cap [S,N], \]
which is equivalent to the assertion that
	\[ R \cap S \subseteq [R,N] \cap [S,N]. \hspace{3mm}\Box\]

\vspace{3mm}

\noindent It is worth noting that if $K$ is $N/RS$-Cockcroft, then $i = 0$,
and so the natural map

\[H_2 \,N/R \oplus H_2 \,N/S \rightarrow H_2 \,N/RS\]

\noindent is injective.  (Compare \cite{gh1,gh2}.)

Taking $N = RS$, the theorem yields a result on asphericity of two-complexes.
All notation is that of the theorem.

\vspace{3mm}

\noindent {\bf Corollary} {\it
$K$  is aspherical if and only if
	\begin{enumerate}
		\item $K_r$ is $RS/R$-Cockcroft,
		\item $K_s$ is $RS/S$-Cockcroft, and
		\item $R \cap S \subseteq [R,R][R,S] \cap [S,S][R,S].\   \Box$
	\end{enumerate}}

\vspace{3mm}

\noindent  This result can also be deduced from  \cite[Theorem 1]{gr}  and an
exact sequence due to R.\ Brown \cite{b}.

\section{Examples}

The condition 3 of the theorem is of interest in the study of the
third integral homology of groups.  The four-term Hopf sequence displayed in
the Introduction shows that elements in third homology are the residues of
homologically trivial spherical maps when the homotopy action of the
fundamental group is trivialized.  It is relatively easy to produce
homologically trivial spherical maps in specific examples.  To decide whether
a spherical homotopy class survives when the fundamental group action is
trivialized is more difficult; this problem involves the internal structure of
$\pi_2$.

A combinatorial approach to this problem was introduced in \cite{bg}.  Among
the results there is a ``relative Hopf formula'' for third homology
\cite[Corollary 5.5]{bg} in the form of an exact sequence

\[ H_3 F/R \oplus H_3 F/S \stackrel{j}{\rightarrow} H_3 F/RS
\stackrel{e}{\rightarrow} \frac{[R,F] \cap [S,F]}{[R,S][F,R \cap S]}
\rightarrow 0.\]

\noindent The value of this sequence begins with the fact that the homology of
one-relator groups is completely understood \cite{l}.  One may think of the
image of $i$ as representing ``obvious'' elements of $H_3 F/RS$.  Results on
the kernel of $j$ appear in \cite{deg}.  The image of $e$ carries the
``nonobvious'' elements of $H_3 F/RS$.  If $R \cap S \subseteq [R,F] \cap
[S,F]$, then every element of $R \cap S$ determines an element in the image of
$e$.  This is relevant to the search for nonobvious elements in $H_3 F/RS$
as it is easier to produce elements of $R \cap S$ than of the subgroup $[R,F]
\cap [S,F]$.

As in \cite[Section 6]{bg}, suitable conditions on $R \cap S$ facilitate the
detection of nonobvious elements of $H_3 F/RS$.  The lower central series of
$F$ provides a good context for calculations.  For a positive integer $n$, let
$F_n$ denote the $n$th term of the lower central series of the free group $F$.
Thus, $F_1 = F$, and $F_{n+1} = [F,F_n]$.  The quotient groups $F_n/F_{n+1}$
are free abelian, with bases determined by the basic commutators of weight
$n$.  (A good reference for this is \cite[p. 149ff]{MHa59}.)  From the
relative Hopf formula, it follows that if $R \cap S \subseteq F_n$, then there
is a homomorphism

\[ e_n : H_3 F/RS \rightarrow \frac{F_n}{[R,S]F_{n+1}}\]

\noindent such that $e_n j = 0$.  The image of $e_n$ therefore detects
nonobvious elements of $H_3 F/RS$.  Elements in the image of $e_n$ are
represented by elements of $[R,F] \cap [S,F]$.

\vspace{3mm}

\noindent {\bf Proposition} {\it Let $({\bf x}:{\bf r},{\bf s})$ be a group presentation and let $F$ be the
free group with basis ${\bf x}$.  Assume that the following conditions hold.
\begin{enumerate}
\item There is a positive integer $n$ such that ${\bf r} \cup {\bf s}
\subseteq F_n$ and
\item the elements $uF_{n+1} \ (u \in {\bf r} \cup {\bf s})$ are linearly
independent in the free abelian group $F_n/F_{n+1}$.
\end{enumerate}
Let $R$ (resp. $S$) denote the normal closure of $\bf r$ (resp. $\bf s$) in
$F$.  Then,
\[R \cap S \subseteq [R,F] \cap [S,F].\]
It follows that $R \cap S \subseteq F_{n+1}$ and that the model of $({\bf
x}:{\bf r},{\bf s})$ is Cockcroft.}

\vspace{3mm}

\noindent {\bf Proof:}  Let  $w \in R \cap S$.  There exist  $r_i \in {\bf
r},\  s_j \in {\bf s},\  w_i,v_j \in F$,  and  $\epsilon _i, \delta_j = \pm 1$
such that
	\[ w = \prod_i w_i r_i^{\epsilon _i} w_i^{-1} = \prod_j v_j
s_j^{\delta _j} v_j^{-1}.  \]
For  $r \in {\bf r}$,  let  
	\[  n_r = \sum_{\textstyle r_i = r} \epsilon _i  \]
and similarly define  $n_s$  for each  $s \in {\bf s}$.  Since
${\bf r} \cup {\bf s} \subseteq F_n$,
	\[  wF_{n+1} = \prod_i r_i^{\epsilon_i}F_{n+1} = \prod_{\textstyle r
\in {\bf r}} r^{\textstyle n_r}F_{n+1}  \]
and hence
	\[  \prod_{\textstyle r \in {\bf r}} r^{\textstyle n_r}F_{n+1} =
\prod_{\textstyle s \in {\bf s}} s^{\textstyle n_s}F_{n+1}.  \]
The linear independence hypothesis now implies that  $n_u = 0$  for all  $u
\in {\bf r} \cup {\bf s}$.  Finally,
	\[  w[R,F] = \prod_i r_i^{\epsilon_i}[R,F] = \prod_{\textstyle r
\in {\bf r}} r^{\textstyle n_r}[R,F] = 1\,[R,F]  \]
so $w \in [R,F]$.  It follows similarly that $w \in [S,F]$.  Since $R
\subseteq F_n$, we also have that $R \cap S \subseteq F_{n+1}$.  Since the
model of a one-relator presentation is Cockcroft, it follows by induction,
compact supports and the theorem that the model $K$ of $({\bf x} : {\bf
r},{\bf s})$ is Cockcroft.  $\Box$

\vspace{3mm}

\noindent This generalizes a result from \cite{d2}, which considers the case
$n=1$.

For the presentations considered in the Proposition, there is the homomorphism
$e_{n+1} : H_3 F/RS \rightarrow F_{n+1}/[R,S]F_{n+2}$.  Examples with $n = 1$
appear in \cite[Examples 6.3,6.4]{bg}.  In each of these examples, the
homomorphism $e_2$ is shown to be nontrivial.

For an example in the case $n=2$, let $F$ be the free group with basis
$\{x,y,z\}$.  Let ${\bf r} = \{[x^a,y]\}$ and ${\bf s} = \{[y^b,z],
[z^c,x]\}$, where $abc \not= 0$.  Since the elements $[x^a,y]F_3 , [y^b,z]F_3
, [z^c,x]F_3$ are linearly independent in $F_2 /F_3 \cong {\bf Z}^{(3)}$, the
Proposition provides that $R \cap S \subseteq [R,F] \cap [S,F] \cap F_3$.
Since $RS \subseteq F_2$, we have the homomorphism

\[e_3 : H_3 F/RS \rightarrow F_3 / [R,S]F_4 = F_3 / F_4 .\]

\noindent Each element $\mu \in R \cap S$ determines an element $\mu F_4$ in
the image of $e_3$.  Now, it can be shown that

\[ \mu = x^a z^c x^{-a} y^b x^a y^{-b} z^{-c}y^b x^{-a} y^{-b} \in R \cap S.\]

\noindent For example, since $x^a$ centralizes all powers of $y$ modulo $R$,

\[ \mu R = x^a z^c x^{-a} x^a y^b y^{-b} z^{-c} y^b y^{-b} x^{-a}R = 1R \]
        
\noindent and so $\mu \in R$.  Similarly, $\mu \in S$.  (It is not too
difficult to produce elements of $R \cap S$ geometrically using equators of
spherical pictures. See \cite{bg,deg} for discussions.)  Manipulations with the
commutator identities yield
\begin{eqnarray*}
\mu F_4 & = & x^a z^c [x^{-a},y^b ] z^{-c} [y^b , x^{-a}]x^{-a}F_4\\
        & = & x^a[z^c,[x^{-a},y^b]]x^{-a}F_4\\
        & = & [[x,y],z]^{abc}F_4 \in F_3 /F_4.
\end{eqnarray*}

\noindent Since $[[x,y],z]F_4$ is
contained in a basis for the free abelian group $F_3 / F_4$, it follows that
$\mu F_4$ is the image under $e_3$ of a nonobvious element of infinite order
in $H_3 F/RS$.  Note that if $K$ is the model of $(x,y,z:{\bf r},{\bf s})$,
then $h: \pi_2 K \rightarrow H_2 K$ is trivial, although both $\pi_2 K$ and
$H_2 K$ are nontrivial.

Take ${\bf r} = \{[x,[x,y^c]]\}$ and ${\bf s} = \{[y,[x,y^c]]\}$ in the
free group $F$ with basis $\{x,y\}$, where $c \not= 0$.  Here, one can show
that
\[\mu = [x,y^c]xy^c x^{-1}[y^c ,x]x[x,y^c ]y^{-c}[y^c ,x]x^{-1} \in R \cap S.\]
The proposition applies with $n=3$, giving the homomorphism \[e_4 : H_3 F/RS
\rightarrow F_4 /F_5.\]  Further, one finds that
\[\mu F_5 = [y,[x[x,y]]]^{c^2}F_5 ,\]
thereby detecting a nonobvious element of infinite order in $H_3 F/RS$.

Huck and Rosebrock \cite{hr} have used the fact that the model of
$(x,y:[x,[x,y^c]],\ [y,[x,y^c]])$ is Cockcroft to show that this presentation
does not have a quadratic isoperimetric inequality.


\begin{thebibliography}{99}


\bibitem{bg} W. A. Bogley and M. A. Guti\'{e}rrez, Mayer-Vietoris sequences
in homotopy of 2-complexes and in homology of groups, J. Pure Appl. Algebra
{\bf 77}  (1992)  39-65.

\bibitem{bd} J. Brandenburg and M. Dyer, On J. H. C. Whitehead's aspherical
question I,  Comment. Math. Helv. {\bf 56}  (1981)  431-446.

\bibitem{bds} J. Brandenburg, M. Dyer and R. Strebel, On J. H. C. Whitehead's
aspherical question II,  in: {\em Low Dimensional Topology}, S. Lomonaco, ed.,
Contemp. Math. {\bf 20} (1983)  65-78.

\bibitem{b}  R. Brown, Coproducts of crossed $P$-modules:  Applications to
second homotopy groups and the homology of groups, Topology {\bf 23} (1984)
337-345. 

\bibitem{c}  W. H. Cockcroft, On two-dimensional aspherical complexes, Proc.
London Math. Soc. {\bf 4} (1954) 375-384.

\bibitem{deg} A. J. Duncan, G. Ellis and N. D. Gilbert, A Mayer-Vietoris
sequence in group homology and the decomposition of relation modules,
preprint, Heriot-Watt University and University College, Galway, 1992.

\bibitem{d1}  M. N. Dyer, Crossed modules and the second homotopy modules of
two-complexes, in: {\em Combinatorial Group Theory and Homotopy in Low
Dimensions}, C. Hog-Angelloni, W. Metzler, and A. J. Sieradski, editors,
London Math. Soc. Lecture Note Series  (Cambridge University Press)  to
appear.

\bibitem{d2}  M. N. Dyer, Cockcroft 2-complexes, preprint, University of
Oregon, 1993.  

\bibitem{g} S. M. Gersten, Dehn functions and $\ell_1$-norms of finite
presentations, in: {\em Algorithms and Classification in Combinatorial Group
Theory}, G. Baumslag and C. F. Miller III, editors, MSRI publications vol. 23
(1991).

\bibitem{gh1} N. D. Gilbert and J. Howie, Threshold subgroups for Cockcroft
2-complexes, preprint, Heriot-Watt University, 1992.

\bibitem{gh2} N. D. Gilbert and J. Howie, Cockcroft properties of graphs of
2-complexes, preprint, Heriot-Watt University, 1992.

\bibitem{gr}  M. A. Guti\'{e}rrez and J. Ratcliffe, On the second homotopy
group, Quart. J. Math. Oxford (2) {\bf 32} (1981) 45-55.

\bibitem{MHa59} M. Hall, Jr., {\em The Theory of Groups}  (Macmillian, 1959).

\bibitem{h} J. Harlander, Minimal Cockcroft subgroups, Glasgow Math. J., to
appear.

\bibitem{ho1} H. Hopf, Fundamentalgruppe und zweite Bettische Gruppe, Comment.
Math. Helv. {\bf 14} (1941) 257-309. 

\bibitem{ho2} H. Hopf, Beitrage zur Homotopietheorie, Comment. Math. Helv.
{\bf 17} (1945) 307-326.

\bibitem{hr} G. Huck and S. Rosebrock, A bicombing that implies a
sub-exponential isoperimetric inequality, Proc. Edinburgh Math. Soc., to
appear.

\bibitem{hu}  J. Huebschmann, Cohomology theory of aspherical groups and of small cancellation groups, J. Pure Appl. Algebra {\bf 14} (1979) 137-143.

\bibitem{l}  R. C. Lyndon, Cohomology theory of groups with a single defining
relation, Ann. of Math. {\bf 52}  (1950)  650-665.

\bibitem{sjp91} S. J. Pride, Identities among relations of groups
presentations, in: {\em Proceedings of the Workshop on Group Theory for a
Geometrical Viewpoint, Trieste 1990}, E.  Ghys, A. Haefliger and A. Verjovsky,
editors  (World Scientific Publishing, Singapore 1991) 687-717.

\bibitem{sjp92} S. J. Pride, Examples of presentations which are minimally
Cockcroft in several defferent ways, preprint, University of Glasgow, 1992.

\end{thebibliography}
\end{document}